\documentclass[letterpaper, 10pt, conference]{ieeeconf}
\IEEEoverridecommandlockouts   
\usepackage{amsmath, amssymb, mathtools, amsthm}
\usepackage{graphicx}
\usepackage{booktabs}
\usepackage{cite}
\usepackage{mathptmx}          
\let\labelindent\relax   
\usepackage{enumitem}
\usepackage{dsfont}
\usepackage{xcolor}
\newif\ifarxiv\arxivtrue   

\theoremstyle{plain}
\newtheorem{theorem}{Theorem}
\newtheorem{proposition}{Proposition}
\newtheorem{lemma}{Lemma}

\theoremstyle{definition}
\newtheorem{assumption}{Assumption}
\newtheorem*{problem}{Problem}
\newtheorem{example}{Example}
\theoremstyle{remark}
\newtheorem{remark}{Remark}

\usepackage{algorithm,algorithmic,float}

\newcommand{\R}{\mathbb{R}}
\newcommand{\Pp}{P^{\perp}}
\newcommand{\Ob}{\mathcal{O}}
\newcommand{\rng}{\mathrm{range}}

\title{\LARGE\bf Sparse Fault Detection in Linear Systems:\\ Equivalence of GLR Tests and $\ell_1$-regularized Least Squares}

\author{Rajasekhar Anguluri$^{1}$%
\thanks{$^{1}$Systems and Control Engineering, Indian Institute of Technology Bombay, Mumbai, India.
        {\tt\small rajanguluri@iitb.ac.in}}}

\begin{document}
\maketitle
\thispagestyle{empty}
\pagestyle{empty}
\begin{abstract}
This paper studies the detection of an impulsive fault in a linear system with unknown initial state and Gaussian measurement noise. The fault incipient time and input channel are unknown. We derive a formula for the fraction of a fault that remains visible after eliminating the initial state. For stable systems, the number of poorly visible fault times stays bounded as the record grows. For diagonalizable systems with undamped modes, a fixed fraction of the record remains poorly visible. We show that the generalized likelihood ratio test detects a fault exactly when a weighted $\ell_1$-regularized least-squares problem (lasso) returns a nonzero estimate. The test threshold sets the lasso penalty and can be calibrated offline. Simulations confirm false-alarm control and show that identifying the faulty channel can be harder than detecting the fault.
\end{abstract}

\vspace{2mm}

\section{Introduction}\label{sec:intro}
Model-based fault and attack detection is a mature subject with renewed relevance in cyber-physical infrastructure, such as power systems and heavy machinery, where an actuator can malfunction or be attacked for a few samples while the plant state is not measured directly. We consider the following basic setup: a linear system is observed through noisy outputs over a finite window; the state at the beginning of the window is unknown; and at an unknown time, through an unknown input channel, the system receives a short impulsive fault. The task is to decide whether a fault occurred, with a guaranteed false-alarm probability that crucially does not depend on the initial state.

\smallskip
\noindent \textbf{Theoretical motivation.} Detectors for unknown input (e.g., faults and attacks) are predominantly chi-square ($\chi^2$) detectors \cite{MWS15,MR17}. These detectors treat the unknown input as an arbitrary vector. Practical faults are rarely arbitrary. Johnson \cite{johnson1998identification} models them as the output of a known linear generator driven by sparse impulses, like step, a ramp, a pulse or a sine whose coefficients jump at a few unknown times. A fault is then a sparse combination of known signatures, one for each onset time and input channel, and the natural detector is a sparse generalized likelihood ratio (GLR) test. 

\smallskip 
We study the simplest sparse structure, a single impulse, and ask what the sparse GLR test gains over the well-known chi-square detector. Answering this question draws on three fields: dynamical systems (fault signatures), linear algebra (projection), and statistical signal processing (invariant tests).

\smallskip
\noindent \textbf{Contributions.} For this fault model we study generalized likelihood ratio (GLR) tests and make three contributions:
\begin{enumerate}
\item \emph{Visibility profile.} The fraction $\rho_a^2$ of a fault signature that survives the data projection has a closed form through observability Gramians (Lemma~\ref{thm:profile}). The hidden fraction is bounded by a geometrically decaying function of $\tau$ for strictly stable systems, and $\rho_a^2$ approaches $(\tau+1)/N$ for diagonalizable systems with all eigenvalues on the unit circle (Theorem~\ref{thm:long}), with fault time $\tau$.
\item \emph{Lasso form.} The invariant GLR detects a fault exactly when a weighted lasso with an unpenalized nuisance block returns a nonzero solution, and the regularization parameter is the threshold. For unknown noise variance the same holds for the square-root lasso (Lemma~\ref{lem:lasso}).
\item \emph{Exact calibration.} The null law depends on the system only through the Gram matrix of the projected atoms and, for unknown variance, the parity-space dimension, and neither on the initial state nor on the variance. Thus the threshold can be calibrated offline (Theorem~\ref{thm:pivot}).
\end{enumerate}
\smallskip

Simulations, including a two-input, two-output example, confirm the predictions. \ifarxiv Appendix~I collects the proofs and Appendix~II additional simulations.\fi

\smallskip
\noindent \textbf{Related literature.} The unknown initial state is a nuisance parameter. In the analytical redundancy literature, this state is removed by projecting the outputs onto the parity space, the orthogonal complement of the range of the observability matrix \cite{CW84,BN93}; Fouladirad and Nikiforov \cite{FN05,FFN08} derive optimal tests, invariant to the initial state, for a fault at a known time. For an unknown time, the generalized likelihood ratio (GLR) of Willsky and Jones \cite{WJ76} maximizes over the onset, and the threshold of the resulting maximum is usually set by simulation. In the security literature, Pasqualetti, D\"orfler and Bullo \cite{PDB13} characterize undetectable attacks as inputs whose output equals the free response for some initial state. Fawzi, Tabuada and Diggavi \cite{FTD14} use $\ell_1$ decoding to handle sparse attacks. Fuchs \cite{Fuchs10} observed that, for a dictionary of atoms in Gaussian noise, the lasso detects a signal exactly when the GLR test does, with the regularization parameter as the threshold. Selective inference \cite{LTTT14,TLT16,TTLT16} gives exact $p$-values for the first variable entering a lasso path.

\smallskip 
Our lasso form in Section~\ref{sec:lasso} builds on \cite{Fuchs10}. We extend it to tests that are invariant to an unknown initial state and noise variance, and we add a closed form for the visible fraction of each fault signature, which shows which fault times are detectable at all.

\section{Problem Setup and Preliminaries}\label{sec:setup}
For a matrix $M$, $\rng(M)$ is its range. Without a subscript, $\|\cdot\|$ is the Euclidean norm of a vector and the spectral norm of a matrix; for a vector $v$, $\|v\|_1$ and $\|v\|_\infty$ are the $\ell_1$ and $\ell_\infty$ norms. For a symmetric matrix, $\lambda_{\min}$ and $\lambda_{\max}$ are its extreme eigenvalues and $M\preceq M'$ is the Loewner order. $\epsilon_j$ is the $j$th unit vector, $\mathbf 1_N$ is the all-ones vector, and $Q(c)=\Pr[z\geq c]$ is the right tail probability of the standard normal $\mathcal N(0,1)$. The \emph{null distribution} of a statistic is its law under $H_0$.

\subsection{System dynamics and attack model}
Consider the standard discrete time linear system: 
\begin{align*}
    x_{k+1}&=Ax_k+Bu_k \qquad  \text{and} \qquad   y_k=Cx_k+e_k, 
\end{align*}
where $x_k\in\R^n$, $u_k\in\R^m$, and $y_k\in\R^p$. The initial state $x_0$ is unknown
and the measurement errors $e_k\sim\mathcal N(0,\sigma^2I_p)$ are i.i.d. The triple $(A,B,C)$ is known and $(A,C)$ observable. The nominal input $u_k$ is zero. 

A fault or an attack is an impulse of unknown amplitude $\mu\neq 0$ that hits
an unknown input channel $j\in\{1,\dots,m\}$ at an unknown time
$\tau\in\{0,\dots,N-2\}$, that is,
\begin{equation*}
  u_k \;=\;
  \begin{cases}
    \mu\,\epsilon_j, & k=\tau,\\[2pt]
    0, & k\neq\tau,
  \end{cases}
\end{equation*}
where $\epsilon_j\in\R^m$ is the $j$th unit vector. We observe $y_0,\dots,y_{N-1}$ with $N\ge n$, so the fault can occur at any
$\tau\in\{0,\dots,N-2\}$. We call the pair $a=(\tau,j)$ a \emph{fault index} and
write $\mathcal A$ for the set of all fault indices, so $|\mathcal A|=m(N-1)$.

Stack the outputs
$y=(y_0,\dots,y_{N-1})\in\R^{Np}$ and errors $e=(e_0,\dots,e_{N-1})\sim\mathcal N(0,\sigma^2I)$ to get the static linear model: 
\begin{equation}\label{eq:model}
  y \;=\; \Ob x_0 \;+\; \mu\, d_a \;+\; e.
\end{equation}
Here $\Ob^\top =\begin{bmatrix}C^\top &(CA)^\top&\ldots& (CA^{N-1})^\top\end{bmatrix}\in\R^{n\times Np}$ is the observability matrix and 
has full column rank for $N\ge n$.
The vector $d_a\in\R^{Np}$ (call as the \emph{signature} of the
fault $a$) is the output response to a unit impulse
$u_\tau=\epsilon_j$:
\[
  (d_a)_k=\begin{cases}CA^{k-\tau-1}b_j, & k>\tau,\\ 0, & k\le \tau,\end{cases}
\]
where $b_j$ is the $j$th column of $B$. 

\subsection{Problem statement}\label{sec:prob_statement}
Since $x_0$ is unknown, $\Ob x_0$ can be any vector in $\rng(\Ob)$. A fault with
signature $d_a$ close to $\rng(\Ob)$ can be explained, almost entirely, by a
different $x_0$. This is the central difficulty of the problem, and the
rest of the paper quantifies this difficulty. We want to decide whether a fault has occurred.

\begin{problem}[Invariant fault detection]\label{prob:main}
Given $y$ in \eqref{eq:model}, test
\[
\begin{aligned}
 \text{(nominal)} \qquad H_0&:\ \mu=0,\quad \\
 \text{(fault)} \qquad  H_1&:\ \mu\neq0 \quad \text{ for some } a\in\mathcal A,
\end{aligned}
\]
without knowing $x_0$ and, in a variant, without knowing $\sigma$. Design a statistical test
such that
\begin{enumerate}[label=(\roman*)]
  \item its false-alarm probability is at most a prescribed value $\alpha$ for every $x_0$
        (and every $\sigma$ in the variant), and
  \item its threshold can be computed from $(A,B,C,N)$ only.
\end{enumerate}
Furthermore, quantify for each $a\in\mathcal A$ the part of the fault that the $x_0$ hides. Requirements (i) and (ii) are addressed in Section~\ref{sec:lasso}. The hidden part of a fault is in
Section~\ref{sec:masking}.
\end{problem}

In the problem above, ``for some $a\in\mathcal A$'' means that, under $H_1$, exactly one of the $|\mathcal A|=m(N-1)$ faults, a pair of input channel and fault time, has occurred. Hence $H_1$ is composite, with one model $y=\Ob x_0+\mu\,d_a+e$ for each index $a$, with unknown amplitude $\mu\neq0$.

\smallskip
\textbf{Performance.} A test $\delta$ declares a fault when $\delta(y)=1$. Its false-alarm probability is $P_{\mathrm{FA}}=\Pr[\delta(y)=1|\,(x_0,\sigma^2), H_0]$, its probability of detection is $P_{\mathrm D}=\Pr[\delta(y)=1|\,(x_0,\sigma^2), H_1]$ for fault parameters $(\mu\ne 0, a)$, and requirement~(i) reads $P_{\mathrm{FA}}\le\alpha$ for every $x_0$ (and $\sigma$).

\section{Generalized Likelihood Ratio Tests}\label{sec:glr}
We consider GLR tests for known and for unknown noise variance $\sigma$. We first eliminate $x_0$ from the data. 

\subsection{Eliminating the initial state}\label{sub:proj}
Let $\Pp$ be the orthogonal projector onto the \emph{parity space} $\rng(\Ob)^\perp$,
which has dimension $r=Np-n$. Since $\Pp \Ob=0$, projecting $y$ \eqref{eq:model} eliminates
$x_0$:
\begin{equation}\label{eq:projected}
    \Pp y=\mu\,\Pp d_a+\Pp e,
\end{equation}
where $\Pp e\sim\mathcal N(0,\sigma^2\Pp)$. Moreover, $\Pp y$ is a maximal invariant\footnote{A test $\phi(y)$ is
\emph{invariant} if $\phi(y+\Ob v)=\phi(y)$ for all $v\in\R^n$, so its decision does not
depend on $x_0$. A test is invariant iff the test is a function of $\Pp y$: two data vectors
with the same projection differ by a vector in $\rng(\Ob)$, so an invariant test cannot
tell them apart.} for the shifts $y\mapsto y+\Ob v$, $v\in\R^n$ \cite{FN05}.

Only the part of a fault signature that no initial state $x_0$ can imitate remains in the projected data $\Pp y$. For each
$a\in\mathcal A$ we define its \emph{visible length} and \emph{visible fraction},
\begin{equation}\label{eq:atoms}
    w_a=\|\Pp d_a\| 
    \quad \text{and }\quad 
    \rho_a^2=\frac{\|\Pp d_a\|^2}{\|d_a\|^2}\in[0,1],
\end{equation}
and call $a\mapsto\rho_a^2$ the \emph{visibility profile}. If $w_a=0$, then
$d_a\in\rng(\Ob)$ and hence the fault is indistinguishable from $x_0$. We exclude this case in
Assumption~\ref{ass:gen}. 

\smallskip 
For $w_a>0$, the \emph{projected atom} $z_a=\Pp d_a/w_a$ is a
unit vector, and the projected data
\begin{equation}\label{eq:reduced}
    \Pp y=(\mu\,w_a)\, z_a+\Pp e.
\end{equation}
Thus, compared with the model with $x_0$ known, the overall fault magnitude is scaled down from $\|d_a\|$ to
$w_a=\rho_a\|d_a\|$. So the problem becomes detecting an unknown multiple of one atom
$z_a$ from the dictionary $\{z_a\}_{a\in\mathcal A}$, with $a$ unknown.

\begin{lemma}[Geometry of visibility]\label{lem:angle}
Let $\theta_a\in[0,\pi/2]$ be the angle between $d_a$ and the subspace $\rng(\Ob)$,
that is, the smallest angle between $d_a$ and a nonzero vector of $\rng(\Ob)$. Then
$\rho_a=\sin\theta_a$.
\end{lemma}
\begin{proof}
Split $d_a=p+\Pp d_a$ with $p=(I-\Pp)d_a\in\rng(\Ob)$ and $\Pp d_a\perp\rng(\Ob)$. For every vector $s\in\rng(\Ob)$ we have that $\langle d_a,s\rangle=\langle p,s\rangle$, so by the
Cauchy--Schwarz inequality $\cos\angle(d_a,s)\le\|p\|/\|d_a\|$, with equality when $s=p$.
Hence $\cos\theta_a=\|p\|/\|d_a\|$. Since $p$ and $\Pp d_a$ are orthogonal and
$d_a=p+\Pp d_a$, the Pythagorean identity gives
\begin{align*}
    \|\Pp d_a\|^2\!=\!\|d_a\|^2\!-\!\|p\|^2\!=\!\|d_a\|^2\big(1\!-\!\cos^2\theta_a\big)\!=\!\|d_a\|^2\sin^2\theta_a.
\end{align*}
Dividing by $\|d_a\|^2$ yields $\rho_a^2=\sin^2\theta_a$, and since $\rho_a\ge0$ and
$\sin\theta_a\ge0$ for $\theta_a\in[0,\pi/2]$, we conclude $\rho_a=\sin\theta_a.$
\end{proof} 

\begin{figure*}[t]
\centering\includegraphics[width=0.95\textwidth]{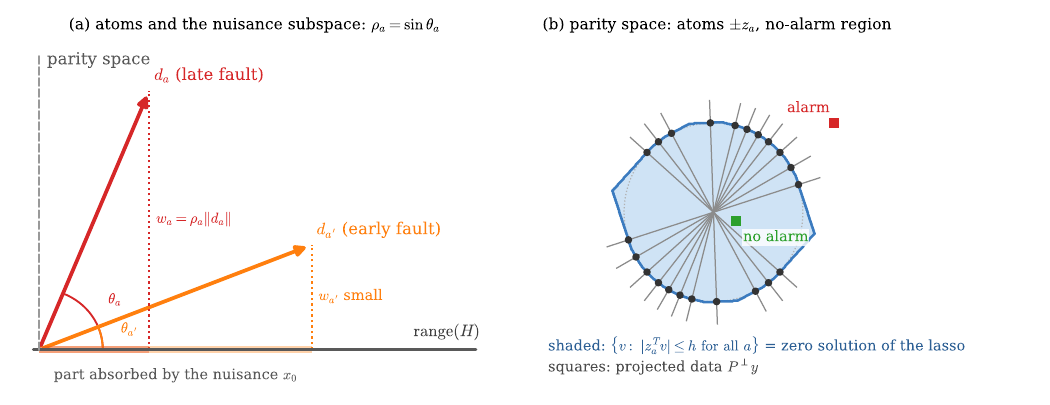}
\caption{Geometry of the detection problem. (a) In $\R^{Np}$ the unknown initial state moves the data freely along $\rng(\Ob)$. Of a fault atom $d_a$ only the component orthogonal to $\rng(\Ob)$ is visible; its length is $w_a=\rho_a\|d_a\|$ with $\rho_a=\sin\theta_a$. (b) In the parity space the normalized visible components are unit vectors $\pm z_a$. The shaded region is the set of projected data for which the lasso of Section~\ref{sec:lasso} returns zero, that is, the no-alarm region of the GLR test with threshold $h$.}\label{fig:geom}
\end{figure*}

\smallskip 
Figure~\ref{fig:geom}(a) illustrates the lemma. The unknown $x_0$ moves the data $y$ only
along $\rng(\Ob)$, so $x_0$ can mimic the part of a fault $d_a\in \rng(\Ob)$. The visible part is the perpendicular component, whose length relative to
$\|d_a\|$ is $\sin\theta_a$. Two corner cases are worth studying. First, if $d_a \in \rng(\Ob)$, a change
in $x_0$ produces exactly the same data, and $\rho_a=0$. Second, if $d_a \in \rng(\Ob)^\perp$, no $x_0$ can mimic it, and $\rho_a=1$. The early fault $d_{a'}$ shown in
Figure~\ref{fig:geom} is close to the former, and the late fault $d_a$ is closer to the latter. 

\begin{assumption}[Distinguishable faults]\label{ass:gen}
For every $a\in\mathcal A$, we impose (i) $w_a> 0$; and (ii) for $a\neq a'$, $|z_a^\top z_{a'}|<1$.
\end{assumption}
Both conditions can be checked from the matrices $(A,B,C)$ and do not rely on the fault parameters $(\mu,\tau)$. 

\smallskip 
Condition (i) states that no fault can be explained by $x_0$ alone. Thus, $d_a\notin\rng(\Ob)$, that is, $\theta_a>0$ in
Fig.~\ref{fig:geom}(a). Condition (i) is the finite-horizon, single-impulse version of an
undetectable input in \cite{PDB13}, and if (i) fails, the data $y$ carry no information about
that fault. For $\tau\ge n-1$ it is equivalent to $d_a\neq0$: if $\Ob x=d_a$, then $CA^kx=0$
for $k\le\tau$, so $x=0$ by observability and hence $d_a=0$. Note that $w_a>0$ does not make a fault easy to detect. A small $\rho_a$ means
only a small fraction of the fault is visible (see Section~\ref{sec:masking}). 

\smallskip 
Condition (ii) states
that no two projected atoms $z_a$ are parallel after projection. Two faults with the same signature cannot be distinguished, and would be merged into one atom; this happens, for
instance, if $b_2$ is a multiple of $b_1$. 

\subsection{The GLR test}\label{sub:glr}
Under $H_1$, with fault index $a$ and amplitude $\mu\neq0$, the density of $y$ is
$p(y| \, (x_0,\sigma^2), H_1)=\mathcal N(\Ob x_0+\mu d_a,\sigma^2I)$, and under $H_0$ the density
$p(y| \, (x_0,\sigma^2), H_0)=\mathcal N(\Ob x_0,\sigma^2I)$. Since $H_1$ is composite and $x_0$ is
unknown, we use the \emph{generalized likelihood
ratio} (GLR) statistic in log form:
\begin{equation}\label{eq:glr}
  \ell(y)=\log\frac{\max_{a\in\mathcal A}\ \sup_{x_0,\mu,\sigma^2}\ p(y\,|\,x_0,\mu,a,\sigma^2)}
                  {\sup_{x_0,\sigma^2}\ p(y\,|\,x_0,0,\sigma^2)},
\end{equation}
where $\sigma^2$ is maximized only in the unknown-variance case, and the supremum over
$\mu\neq0$ equals the supremum over $\mu\in\R$ because $p$ is continuous in $\mu$. The GLR
test is
\begin{equation}\label{eq:glrt}
  \delta(y)=\begin{cases}1\ \text{(alarm)}, & \ell(y)\ge \eta,\\ 0, & \ell(y)<\eta.\end{cases}
\end{equation}
Here $\eta\geq 0$ is such that $P_{\mathrm{FA}}=\mathrm{Pr}[\ell(y)\geq \eta|\,(x_0, \sigma),H_0]\le\alpha$. Because $H_1$ is composite, optimal tests usually do not exist, and the GLR test is sub-optimal.
To evaluate $\ell(y)$ we need the residuals of the least-squares fits under each hypothesis.

\begin{lemma}[Residuals]\label{lem:glr}
For each fault index $a\in\mathcal{A}$,
\begin{align*}
    \min_{x_0,\mu}\|y-\Ob x_0-\mu d_a\|^2
    &=\|\Pp y\|^2-(z_a^\top y)^2,\\
    \min_{x_0}\|y-\Ob x_0\|^2
    &=\|\Pp y\|^2.
\end{align*}
\end{lemma}

\begin{proof}
Consider the first minimization problem. For a fixed $\mu$, least-squares minimization
over $x_0$ gives
\[
    \min_{x_0}\|y-\Ob x_0-\mu d_a\|^2
    =\|\Pp(y-\mu d_a)\|^2
    =\|\Pp y-\mu w_a z_a\|^2, 
\]
where the final equality follows from \eqref{eq:reduced}. Because $z_a$ is a unit vector and $w_a>0$, the amplitude $\mu^\star=(z_a^\top\Pp y)/{w_a}$ minimizes $\|\Pp y-\mu w_a z_a\|^2$. From these results we get 
\[
    \min_{x_0,\mu}\|y-\Ob x_0-\mu d_a\|^2
    =\|\Pp y\|^2-(z_a^\top\Pp y)^2.
\]
Because $z_a\in\rng(\Ob)^\perp$, we have
$z_a^\top\Pp y=z_a^\top y$, proving the first identity.
For the fault-free model, setting $\mu=0$ in the
least-squares problem above gives the second identity.
\end{proof}

\smallskip
We evaluate \eqref{eq:glr}. Let $\ell_{\mathrm{proj}}(y)$ be defined as in \eqref{eq:glr}
with the projected data $\Pp y$ in place of $y$. The densities of $\Pp y$ under $H_0$ and $H_1$ are Gaussian and do not depend on $x_0$. The next proposition evaluates $\ell_{\mathrm{proj}}$ and $\ell$, and shows that projecting the data loses no information for the GLR test. 

\begin{proposition}[The GLR test]\label{prop:glr}
Let Assumption~\ref{ass:gen}(i) hold, and let $r=Np-n\ge2$. Define the statistics
\begin{equation}\label{eq:T}
  T=\max_{a\in\mathcal A}\frac{(z_a^\top y)^2}{\sigma^2},
  \qquad
  U=\max_{a\in\mathcal A}\frac{|z_a^\top y|}{\|\Pp y\|}\in[0,1].
\end{equation}
\begin{enumerate}[label=(\roman*)]
  \item (Known $\sigma$): $\ell=\ell_{\mathrm{proj}}=T/2$.
  \item (Unknown $\sigma$): $\ell_{\mathrm{proj}}=-\tfrac r2\log(1-U^2)$
        and $\ell=\tfrac{Np}{r}\,\ell_{\mathrm{proj}}$.
\end{enumerate}
\end{proposition}

\smallskip
By Proposition~\ref{prop:glr}, $\ell_{\mathrm{proj}}(y)$ is an increasing function of $T$ when $\sigma$ is known and of $U$ when $\sigma$ is unknown. Hence the GLR test of the projected data, $\ell_{\mathrm{proj}}(y)\ge\eta$, is
\begin{equation}\label{eq:tests}
  \delta(\Pp y)=
  \begin{cases}
    \mathds 1\{T\ge\gamma\}, & \sigma\text{ known},\ \ \gamma=2\eta,\\[2pt]
    \mathds 1\{U\ge u\}, & \sigma\text{ unknown},\ \ u=\sqrt{1-e^{-2\eta/r}}.
  \end{cases}
\end{equation}
Proposition~\ref{prop:glr} shows that the log-likelihood ratios of $y$ and of $\Pp y$ are the same for known $\sigma$, and for unknown $\sigma$ they are equal up to the scaling factor $Np/r$. Hence the GLR test statistics are the same $T$ and $U$ for the unprojected data, and the only change is in the thresholds determined by $\eta$ in \eqref{eq:glrt}. In practice we choose $\gamma$ and $u$ directly as the smallest values with $\Pr[T\ge\gamma]\le\alpha$ and $\Pr[U\ge u]\le\alpha$ under $H_0$ (Section~\ref{sec:lasso}). See Remark~\ref{rem:stats} for how $U$ relates to $T$.


\begin{example}[Integrator]\label{ex:integrator}
Let $A=B=C=1$, so $\Ob=\mathbf 1_N$ and $d_\tau=(0,\dots,0,1,\dots,1)^\top$ with
$q=N-\tau-1$ ones. Set $t=\tau+1$. Then $\Pp=I-\mathbf 1_N\mathbf 1_N^\top/N$,
$w_\tau^2=\|d_\tau\|^2-q^2/N=tq/N$, hence $\rho_\tau^2=t/N$ and $z_\tau^\top y=\sqrt{tq/N}\,(\bar y_{>\tau}-\bar y_{\le\tau})$, where $\bar y_{>\tau}$ and $\bar y_{\le\tau}$ are the sample means after and up to time
$\tau$. Under $H_0$, $z_\tau^\top y/\sigma$ is the standardized difference of two sample
means, so $T$ is the classical scan statistic for a change in mean \cite{BN93}. An early
fault is poorly visible because only $t$ samples remain to estimate the unknown pre-fault level.
\end{example}

\begin{remark}[The two statistics]\label{rem:stats}
Under $H_0$, $\hat\sigma^2=\|\Pp y\|^2/r$ is unbiased for $\sigma^2$, and $rU^2=\max_a(z_a^\top y)^2/\hat\sigma^2$. So $rU^2$ is $T$ with $\sigma^2$ replaced by $\hat\sigma^2$, and $U$ is the cosine of the angle between $\Pp y$ and the closest atom $\pm z_a$. $U$ is scale invariant, so its null distribution does not depend on $\sigma$ (Theorem~\ref{thm:pivot}).
\end{remark}

\begin{remark}[Comparison with $\chi^2$ fault/attack detectors]\label{rem:chi2}
In the well-known $\chi^2_r$ attack detector for cyber-physical systems \cite{MWS15,MR17}, the test statistic is $\|\Pp y\|^2/\sigma^2$. This is the value of $2\ell$ in \eqref{eq:glr} if the fault/attack is allowed to be any vector of the parity space $\rng(\Ob)^\perp$ instead of an impulse fault (our case). Here $r=Np-n$ is the number of degrees of freedom, the dimension of the parity space $\rng(\Ob)^\perp$. 

\smallskip 
Our fault/attack model is an impulse, a single atom of that space. The $\chi^2_r$ detector therefore spreads an impulse over all $r$ dimensions and cannot return $(\tau,j)$, whereas $T$ matches one atom and returns the maximizing pair. With an arbitrary mean, $\sigma$ is not identifiable from $\Pp y$, so the $\chi^2_r$ detector needs $\sigma$ given, whereas the impulse structure lets $U$ handle unknown $\sigma$. Detectors built on Kalman innovations assume a prior on $x_0$, and their false-alarm probability is guaranteed only when $x_0$ follows that prior; ours holds for every $x_0$. \qed 
\end{remark}

\section{Structure and Calibration of the Test}\label{sec:lasso}
The test statistics $T$ and $U$ depend on the model \eqref{eq:model} only through the unit atoms $z_a$, so the results in this section hold for any collection of atoms, that is, a finite frame \cite{SC17}. The system dynamics enters only in Section~\ref{sec:masking}. We show that the correlations among the atoms, collected in the Gram matrix $G=[z_a^\top z_{a'}]$, determine the geometry of the test as a lasso, the null distribution of its statistics, and hence the thresholds.

\subsection{Lasso form}\label{sub:lasso}
Let $D=[d_a]_{a\in\mathcal A}\in\R^{Np\times|\mathcal A|}$ be the matrix whose columns are the fault signatures in \eqref{eq:model}. For $h,\lambda>0$, consider the lasso problems
\begin{align}
&\textbf{(P1)}\ \ \underset{x_0,\,\psi}{\arg\min}\ \tfrac12\|y-\Ob x_0-D\psi\|_2^2+h\sum_{a}w_a|\psi_a|,\label{eq:wlasso}\\
&\textbf{(P2)}\ \ \underset{x_0,\,\psi}{\arg\min}\ \|y-\Ob x_0-D\psi\|_2+\lambda\sum_{a}w_a|\psi_a|.\label{eq:sqrt}
\end{align}
Let $(\hat x_0,\hat\psi)$ denote a solution of (P1) or (P2). Here $D\psi=\sum_{a}\psi_a d_a$, so $\psi_a$ is the weight on the fault $d_a$. The variable $x_0$ is not penalized, and $w_a$ is the visible length in Fig.~\ref{fig:geom}(a). Our model \eqref{eq:model}, with a single fault $a$ of amplitude $\mu$, is the special case $\psi=\mu\,\mathbf e_a$, where $\mathbf e_a$ is the $a$-th standard basis vector. A solution with $\hat\psi\neq0$ declares a fault. Lemma~\ref{lem:lasso} below shows that this reproduces the GLR tests of Proposition~\ref{prop:glr}.

\begin{lemma}\label{lem:lasso}
Under Assumption~\ref{ass:gen}, $\hat\psi=0$ is a solution of (P1) if and only if $h\ge\max_a|z_a^\top y|$, that is, $T\le h^2/\sigma^2$. If $\Pp y\neq0$, $\hat\psi=0$ is a solution of (P2) if and only if $\lambda\ge U$.
\end{lemma}

Hence the GLR tests in \eqref{eq:tests} reject $H_0$ exactly when $\hat\psi\neq0$, with $h=\sigma\sqrt\gamma$ in (P1) and $\lambda=u$ in (P2). Geometrically (see Fig.~\ref{fig:geom}(b)), $z_a^\top y=z_a^\top\Pp y$, so no fault is detected exactly when $\Pp y$ lies in the set $\{v\in\rng(\Ob)^\perp:|z_a^\top v|\le h\ \forall a\}$. This set is a polytope when the atoms span $\rng(\Ob)^\perp$.

\smallskip 
(\emph{Insight and limits}) Lemma \ref{lem:lasso} gives the penalty term in \textbf{(P1)} and \textbf{(P2)} a false-alarm interpretation. Indeed with $h=\sigma\sqrt\gamma$ or $\lambda=u$ calibrated as in Section~\ref{sub:null}, the lasso has false-alarm probability at most $\alpha$ for every $x_0$, in place of cross-validation. The free block $x_0$ makes the test invariant to the initial state. The weights $w_a$ set the level of each atom: $\hat\psi=0$ is a solution of (P1) if and only if $h\ge\max_aw_a|z_a^\top y|$, so the rule compares $|z_a^\top y|$ with $h/w_a$. With unit weights, all atoms share the level $h$. The lemma does not make the test cheaper, since the maximum is one pass over the atoms whereas a lasso solver iterates. The gain is that the lasso also returns $\hat\psi$ and brings standard solvers and theory.

\subsection{Exact threshold from the Gram matrix}\label{sub:null}
We show that under $H_0$ the laws of $T$ and $U$ are determined by the Gram matrix $G=[z_a^\top z_{a'}]_{a,a'\in\mathcal A}$ of the atoms, and that the law of $U$ does not depend on $x_0$ or $\sigma$, so the thresholds of the GLR tests can be computed once offline.

\begin{theorem}[Null distribution]\label{thm:pivot}
Let $V\in\R^{Np\times r}$ be an orthonormal basis of $\rng(\Ob)^\perp$ and define the vector $z_a=V\zeta_a$. Let $g\sim\mathcal N(0,I_r)$. Under $H_0$, for every $x_0$ and $\sigma$,
\[
T\overset{d}{=}\max_a(\zeta_a^\top g)^2,\qquad \text{and} \qquad U\overset{d}{=}\max_a|\zeta_a^\top v|,
\]
where $v$ is uniform on the unit sphere of $\R^r$.
\end{theorem}

\begin{proof}
Under the hypothesis $H_0$, $\Pp y=\Pp e$ because $\Pp\Ob=0$. Put $g=V^\top e/\sigma\sim\mathcal N(0,I_r)$. Then $z_a^\top y=z_a^\top\Pp e=\sigma\zeta_a^\top g$, which gives the statement for $T$. Also $\|\Pp y\|=\sigma\|g\|$, so $U=\max_a|\zeta_a^\top g|/\|g\|$. A standard Gaussian is rotation invariant, so $v=g/\|g\|$ is uniform on the sphere.
\end{proof}

\smallskip 
By Theorem~\ref{thm:pivot}, thresholds $\gamma$ for $T$ (known $\sigma$) and $u$ for $U$ with $\Pr[T\ge\gamma]\le\alpha$ and $\Pr[U\ge u]\le\alpha$ under $H_0$ give false-alarm probability at most $\alpha$ for all $x_0$ and $\sigma$. For a Gaussian mean test \cite{anguluri2016network} or a $\chi^2$ test, the threshold is an inverse of the $Q$ function or the $\chi^2$ cdf. Here $T$ and $U$ are maxima of correlated variables with no closed-form cdf, so Algorithm~\ref{alg:thr} computes the thresholds offline by Monte Carlo sampling.
\begin{algorithm}[H]
\caption{Offline thresholds $\gamma$ and $u$}\label{alg:thr}
\begin{algorithmic}[1]
\REQUIRE atoms $\zeta_a$, $a\in\mathcal A$; level $\alpha$; number of samples $M$
\FOR{$k=1,\dots,M$}
\STATE draw $g_k\sim\mathcal N(0,I_r)$
\STATE $T_k\leftarrow\max_a(\zeta_a^\top g_k)^2$, \quad $U_k\leftarrow\max_a|\zeta_a^\top g_k|/\|g_k\|$
\ENDFOR
\STATE $\gamma\leftarrow$ value above which a fraction $\alpha$ of $T_1,\dots,T_M$ lie
\STATE $u\leftarrow$  value above which a fraction $\alpha$ of $U_1,\dots,U_M$ lie
\ENSURE $\gamma$ for $T$ (known $\sigma$) and $u$ for $U$
\end{algorithmic}
\end{algorithm}

\smallskip 
\begin{remark}[Closed-form test thresholds]\label{rem:closed}
Sampling can be avoided with closed-form thresholds $c$ with the false alarm probability $\Pr[\max_a|z_a^\top y|/\sigma\ge c]\le\alpha$ under $H_0$, which give $\gamma=c^2$. This test declares a fault when at least one of the $|\mathcal A|$ single-atom tests $|z_a^\top y|/\sigma\ge c$ fires, and each single-atom test has false-alarm probability $2Q(c)$. 

\smallskip 
\emph{Bonferroni} correction \cite{Dunn61} adds these probabilities over the atoms and solves for $2|\mathcal A|Q(c)=\alpha$. The \emph{\v{S}id\'ak} bound in \cite{Sidak67}, the independence approximation used by Fuchs \cite{Fuchs10}, multiplies them as if the atoms were independent and solves for $1-(1-2Q(c))^{|\mathcal A|}=\alpha$. Finally, the \emph{Rice-type tube bound} \cite{Sun93} uses the correlation between neighboring atoms for one input channel and solves $2Q(c)+(\kappa_0/\pi)e^{-c^2/2}=\alpha$, where $\kappa_0=\sum_\tau\arccos(z_\tau^\top z_{\tau+1})$ is the length of the polygonal curve through consecutive atoms on the sphere. 

\smallskip 
All three are conservative when the atoms are correlated. \qed 
\end{remark}

\section{Probability of Detection}\label{sec:masking}
The thresholds of Section~\ref{sec:lasso} fix the false-alarm probability. We now ask how often a fault is detected. Section~\ref{sub:interf} shows that the unknown initial state scales the signal-to-noise ratio by the visible fraction $\rho_a$, and Section~\ref{sub:profile} computes $\rho_a$ as a function of the fault time and the record length.

\subsection{Interference from the initial state}\label{sub:interf}
Projecting out the unknown initial state $x_0$ loses nothing for the GLR test (Proposition~\ref{prop:glr}), but $x_0$ interferes with the fault: the component of $d_a$ in $\rng(\Ob)$ is indistinguishable from $\Ob x_0$ and is removed, so by \eqref{eq:projected} only $\Pp d_a$ remains. A fault of amplitude $\mu$ has visible length $\mu w_a=\mu\rho_a\|d_a\|$, so the visible fraction $\rho_a$ of Section~\ref{sec:glr} is the fraction of the fault that survives, a geometric quantity by Lemma~\ref{lem:angle}. The next proposition states the consequence for the probability of detection $P_{\mathrm D}(a,\mu)$, defined in Section~\ref{sec:prob_statement}.

\begin{proposition}[Interference and probability of detection]\label{prop:pd}
Let $\sigma$ be known and let $\delta(\Pp y)=\mathds 1\{\max_a|z_a^\top y|/\sigma\ge c\}$. For every $x_0$, every fault $a$ and every amplitude $\mu$,
\[
  P_{\mathrm D}(a,\mu)\ \ge\ Q(c-s_a)+Q(c+s_a).
\]
where $Q(\cdot)$ is right tail probability of a standard Gaussian and $s_a=(\mu w_a)/{\sigma}=\rho_a\,(\mu\|d_a\|)/{\sigma}$. 
\end{proposition}

\begin{proof}
Since $z_a\in\rng(\Ob)^\perp$ and $\Pp d_a=w_az_a$, we have $z_a^\top\Ob x_0=0$ and $z_a^\top d_a=w_a$, so $z_a^\top y=\mu w_a+z_a^\top e$, which does not depend on $x_0$. As $z_a$ is a unit vector, $z_a^\top y/\sigma\sim\mathcal N(s_a,1)$. The event $\{\max_{a'}|z_{a'}^\top y|/\sigma\ge c\}$ contains $\{|z_a^\top y|/\sigma\ge c\}$, whose probability is $Q(c-s_a)+Q(c+s_a)$, the two terms being the upper and lower Gaussian tails.
\end{proof}

\smallskip 
The raw signal-to-noise ratio $\mu\|d_a\|/\sigma$ is the fault length in units of the noise standard deviation, which the test would see if $x_0$ were known. Since $x_0$ is unknown, the $T$ test sees the effective ratio $s_a=\rho_a\mu\|d_a\|/\sigma$. A poorly visible fault is thus detected less often. Reaching the same probability of detection needs an amplitude $1/\rho_a$ times larger (discussed in Fig.~\ref{fig:power}B). The bound is conservative, since atoms other than $a$ can also trigger detection.

\subsection{Visibility profile}\label{sub:profile}
In this section we discuss the impact of system dynamics on the signal-to-noise ratio via $\rho_a^2$. First, we compute $\rho_a^2$ from $(A,B,C)$. Let $W_i=\sum_{k=0}^{i-1}(A^k)^\top C^\top CA^k$ be the observability Gramian over $i$ samples. Then $\Ob^\top \Ob=W_N\succ0$ because $N\geq n$ and $(A,C)$ is observable; and $W_i\preceq W_{i'}$ for $i\le i'$. 


\begin{lemma}[Visibility profile]\label{thm:profile}
For a fault index pair $a=(\tau,j)$, with $a \in \mathcal{A}$, let $t=\tau+1$ and $q=N-t$. Then
\begin{equation}\label{eq:profile}
1-\rho_a^2=\frac{b_j^\top W_q\,A^{t}\,W_N^{-1}(A^{t})^\top W_q\,b_j}{b_j^\top W_q\,b_j}.
\end{equation}
\end{lemma}
The proof is a short computation, and we sketch only the high-level details. Equation~\eqref{eq:profile} follows from $d_a=\Ob x_0^\star+\Pp d_a$, where $x_0^\star=W_N^{-1}\Ob^\top d_a$ is the initial state whose fault-free output $\Ob x_0^\star$ is closest to $d_a$. As a result $\Ob x_0^\star$ is the part of $d_a$ that some initial state can reproduce, and $\Pp d_a$ is the residual. By orthogonality of the two parts, we have $1-\rho_a^2=\|\Ob x_0^\star\|^2/\|d_a\|^2$, and substituting $\Ob^\top\Ob=W_N$ and $\Ob^\top d_a=\sum_{i=0}^{q-1}(A^{t+i})^\top C^\top CA^ib_j=(A^t)^\top W_q b_j$ gives \eqref{eq:profile}.

\smallskip 
The role of the fault time $\tau$ in $\rho^2_a$ is apparent from the definition of $d_a$. For $k\ge t$ we have $(d_a)_k=CA^{k-t}b_j$, while $(\Ob x_0)_k=CA^kx_0$, and $(d_a)_k=0$ for $k<t$. If the fault entered at the first sample ($t=0$), $d_a$ would equal $\Ob b_j$, so $d_a$ would be an initial-state change and hence hidden. For $t>0$, $d_a$ vanishes on the first $t$ samples, whereas $\Ob x_0$ vanishes for special $x_0$. A larger $t$ gives more samples of $y$ on which the two must differ, so they are harder to align and $\rho_a^2$ is larger.

\smallskip 
\begin{figure}[t]
\centering\includegraphics[width=\columnwidth]{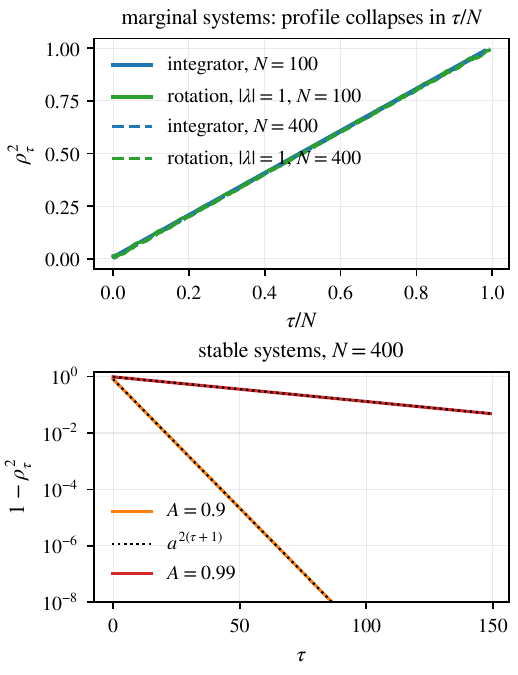}
\caption{Visibility profile $\rho_\tau^2$. Top: marginally stable systems, the integrator and the rotation by $30^\circ$ per step with $B=(1,0)^\top$, $C=(1,0)$; $N=100$ (solid) and $N=400$ (dashed). Bottom: strictly stable scalar systems $A\in\{0.9,0.99\}$ with $B=C=1$, $N=400$; dotted, $a_0^{2(\tau+1)}$.}\label{fig:profile}
\end{figure}

\emph{(Scalar case).} For $n=1$ and $A=a_0\neq0$, the formula \eqref{eq:profile} reduces to $\rho_\tau^2=(1-a_0^{2(\tau+1)})/(1-a_0^{2N})$, independent of $C$ and $b_j$. For the integrator ($a_0=1$) $\rho_\tau$ equals $(\tau+1)/N$, as in Example~\ref{ex:integrator}. For $|a_0|<1$ and large $N$, $\rho_\tau^2$ tends to $1-a_0^{2(\tau+1)}$, so the hidden fraction decays geometrically in $\tau$.

The next result describes the profile for long records. Call an atom \emph{poorly visible} if $\rho_a^2<\rho_0^2$ for a fixed $\rho_0\in(0,1)$.

\begin{theorem}[Visibility for long records]\label{thm:long}\phantom{raj}
\begin{enumerate}[label=(\roman*),leftmargin=1.8em,labelsep=0.4em,itemsep=2pt,topsep=2pt,parsep=0pt]
    \item Let $A$ be Schur stable, $\|A^t\|\le M\lambda^t$ with $\lambda<1$. Let $W_\infty=A^\top W_\infty A+C^\top C$. Then for every $N\ge n$ and every atom,
$1-\rho_a^2\le\kappa_\infty M^2\lambda^{2(\tau+1)}$, with $\kappa_\infty=\lambda_{\max}(W_\infty)/\lambda_{\min}(W_n)$.

\smallskip  
 \item Let $A$ be diagonalizable with all eigenvalues on the unit circle, $(A,C)$ observable, and $b_j\neq0$ for all $j$. Then
\begin{equation}\label{eq:marginal}
\max_{\tau}\ \Big|\rho_{(\tau,j)}^2-\frac{\tau+1}{N}\Big|\ \longrightarrow\ 0\qquad(N\to\infty).
\end{equation}
\end{enumerate}
\end{theorem}

Theorem~\ref{thm:long} describes how the visible fraction ratio $\rho_a^2$ of a single fault $a=(\tau,j)$ depends on the fault time $\tau$, on any input channel $j$, in two regimes. In (i), $\rho^2_a$ tends to one geometrically in the fault time, at a rate $\lambda^{2(\tau+1)}$. Only the first few atoms per input channel are poorly visible, and their number does not grow with $N$ (see Remark~\ref{rem:long}), so a longer record hides no more faults. For $n=1$ the exact hidden fraction $\rho_a^2=a_0^{2t}$, and the poorly visible window has length $\ln(1/(1-\rho_0^2))/(2\ln(1/|a_0|))$. 

\smallskip 
In (ii), system modes neither decay nor grow, and $\rho_a^2\approx(\tau+1)/N$ for every $(A,B,C)$ in the class, as for the integrator of Example~\ref{ex:integrator}. A fault in the first $\rho_0^2N$ samples is therefore poorly visible, however long the record, and the window of poorly visible faults grows in proportion to $N$. By Proposition~\ref{prop:pd}, the probability of detection of such a fault is limited by $s_a=\mu w_a/\sigma$, so for marginally stable plants a longer record helps late faults but not early ones. 

\smallskip 
Fig.~\ref{fig:profile} plots $\rho_a^2$, computed from its definition, and shows the two regimes of Theorem~\ref{thm:long}. For systems with all eigenvalues on the unit circle, $\rho_a^2$ depends on $\tau/N$; for strictly stable systems, the poorly visible window has a fixed length.

\smallskip
\begin{remark}[Limit and count]\label{rem:long}
In Theorem~\ref{thm:long} (i), $1-\rho^2_{(t-1,j)}\to\gamma_{t,j}$ as $N\to\infty$, since $W_q,W_N\to W_\infty$ in \eqref{eq:profile}, where $\gamma_{t,j}=(b_j^\top W_\infty A^tW_\infty^{-1}(A^t)^\top W_\infty b_j)/(b_j^\top W_\infty b_j)$. Moreover, for any $N$, fewer than $\ln\!\big(\kappa_\infty M^2/(1-\rho_0^2)\big)\big/\big(2\ln(1/\lambda)\big)$ atoms per channel are poorly visible. In case~(ii), the poorly visible atoms per channel is $\rho_0^2N+o(N)$ for every such $(A,B,C)$.
\end{remark}

\section{Numerical Studies}\label{sec:sims}
We validate two theoretical results: (i) the exact thresholds give a false-alarm probability $\alpha$ regardless of the initial state $x_0$ (Theorem~\ref{thm:pivot}), and (ii) a fault is seen through its visible fraction $\rho_a$ (Proposition~\ref{prop:pd}). In the last study, we consider fault channel identification. We use $\alpha=0.05$, $N=100$ and the integrator $A=B=C=1$ of Example~\ref{ex:integrator} unless stated, take $\sigma=1$ and estimate probabilities by Monte Carlo runs. 

\smallskip 
The \emph{max test} declares a fault when $T\ge\gamma$, with $\sigma$ known, and the \emph{$U$ test} when $U\ge u$, with $\sigma$ unknown; $\gamma\approx9.4$ and $u\approx0.30$ come from Algorithm~\ref{alg:thr}. We further consider (i) the \emph{oracle test} that knows $x_0$ and uses $\max_a|d_a^\top(y-\mathcal Ox_0)|/\|d_a\|$ on the unprojected data, with its exact threshold; and (ii) the \emph{naive test}, the oracle with $x_0$ wrongly assumed to be $0$.

\begin{figure}[t]
\centering\includegraphics[width=\columnwidth]{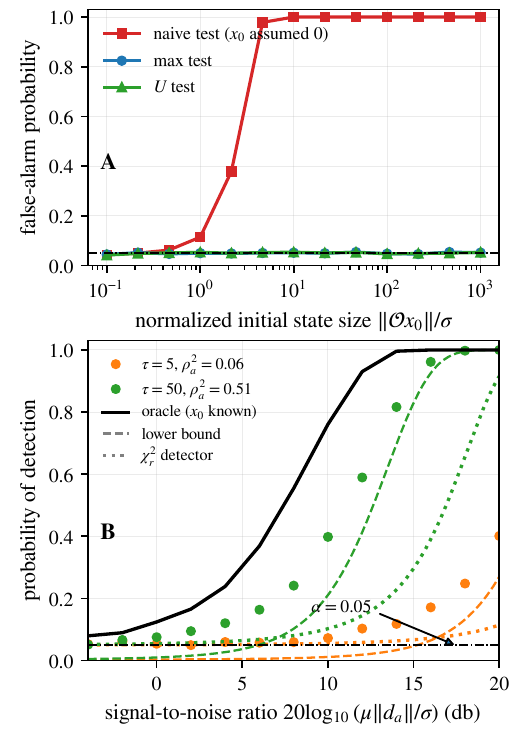}
\caption{Integrator, $N=100$, $\alpha=0.05$, $\sigma=1$. (A) False-alarm probability against the normalized initial state size. (B) Probability of detection for a fault at $\tau=5$ (orange) and $\tau=50$ (green): the max test (dots), its lower bound from Proposition~\ref{prop:pd} (dashed) and the $\chi^2$ detector (dotted). The black curve is the oracle test for $\tau=50$, and the dash-dotted line is $\alpha$.}\label{fig:power}
\end{figure}

\medskip 
\noindent\emph{(False-alarm control)} Fig.~\ref{fig:power}A plots the false-alarm probability against the normalized initial state size $\|\Ob x_0\|/\sigma$, from $0.1$ to $10^3$, with thresholds fixed. Each point is the fraction of $4000$ fault-free trials in which the test detects a fault. The max and $U$ tests stay at $\alpha=0.05$ at every size, as Theorem~\ref{thm:pivot} predicts. The naive test, which assumes $x_0=0$, is near $\alpha$ for small $x_0$ but detects a fault in nearly every trial once $\|\Ob x_0\|/\sigma\gtrsim5$, because it mistakes $x_0$ for a fault. Projecting out $x_0$ avoids this failure.

\smallskip
\noindent\emph{(Price of the unknown $x_0$)} In Fig.~\ref{fig:power}B, $\sigma$ is known and ($x_0$, $\tau$, $\mu$) are unknowns. We run two experiments with fault times at $\tau=5$ and $\tau=50$. The probability of detection of the max test is plotted against the raw signal-to-noise ratio $\mu\|d_a\|/\sigma$ in db ($3000$ trials per point). The test sees the fault at the effective ratio $s_a=\rho_a\mu\|d_a\|/\sigma$ (Proposition~\ref{prop:pd}), since the projection removes the hidden part of the fault. The oracle test, which knows $x_0$, has no such loss. Relative to this test, the effective ratio is lowered by $-10\log_{10}\rho_a^2$ db, where $\rho_a^2=(\tau+1)/N$ for the integrator: $12.2$ db at $\tau=5$ and $2.9$ db at $\tau=50$.

\smallskip 
The dashed curves are the lower bound of Proposition~\ref{prop:pd}; the bound lies below the simulated probabilities because it ignores the other atoms. The oracle curve, drawn for $\tau=50$, differs by at most $0.04$ between fault times. At $\tau=50$, the $\chi^2$ detector discussed in Remark~\ref{rem:chi2} needs about $5$ db more signal than the max test, because it spreads the impulse over $r=99$ directions.

\smallskip
\noindent\emph{(Identifying faults)} Beyond detection, the maximizing atom $(\hat\tau,\hat{j})$, a (time, channel) pair, estimates when the fault happened and which input channel it hit. We ask how often this estimate is correct. The system is a chain of two unit masses with stiffness matrix $K=\left[\begin{smallmatrix}2&-1\\-1&1\end{smallmatrix}\right]$ and damping $\beta K$, $\beta=0.3$. With $\dot x=A_cx+B_cu$ and $y=Cx$, where $A_c=\left[\begin{smallmatrix}0&I\\-K&-\beta K\end{smallmatrix}\right]$, $B_c=\left[\begin{smallmatrix}0\\I\end{smallmatrix}\right]$, $C=[\,I\;\;0\,]$. Input channels ($j=1,2$) are either mass and outputs are positions. We discretize by zero-order hold with step $0.5$ to get Schur stable $(A,B)$. We set $N=200$ and inject impulse fault on channel $1$ at time $\tau$ with $\mu\|d_a\|/\sigma=5$. 

\begin{figure}[t]
\centering\includegraphics[width=\columnwidth]{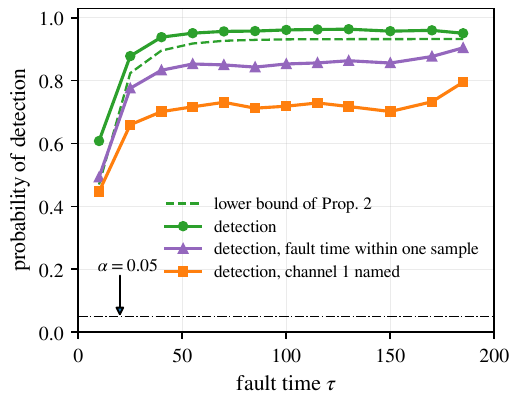}
\caption{Two-mass chain, $\beta=0.3$, $N=200$, $\mu\|d_a\|/\sigma=5$, fault on channel $1$ at time $\tau$, $3000$ trials per $\tau$. Probability that the max test detects the fault (green, with the lower bound of Proposition~\ref{prop:pd} dashed), that it detects and estimates the fault time within one sample (purple), and that it detects and names channel $1$ (orange). The dash-dotted line is $\alpha$.}\label{fig:mimo}
\end{figure}

\smallskip
Fig.~\ref{fig:mimo} plots the probability of detecting the fault with the max test ($\gamma\approx12.3$, $3000$ trails per point) under three criteria. A trial counts as a detection if (i) $T\ge\gamma$ (test detects fault), (ii) $T\ge\gamma$ and the time estimate satisfies $|\hat\tau-\tau|\le1$, and (iii) $T\ge\gamma$ and the channel estimate is $\hat{j}=1$.

\smallskip 
The probability of detection improves with $\tau$ because more of the fault is visible (Theorem~\ref{thm:long}), and levels off at $0.96$, above the bound of Proposition~\ref{prop:pd}. Identifying the channel is harder than detecting the fault or estimating its time, at every $\tau$. At $\tau=55$ the test detects the fault in $95\%$ of trials, with the time within one sample in $85\%$ and the correct channel in only $72\%$. A likely reason is that the projected atoms of the two channels have correlation $0.95$, so the maximum often lands on the other channel.

\section{Conclusion}\label{sec:concl}
We studied the detection of an impulsive fault in the inputs of a linear system with unknown initial state $x_0$, when neither the fault time nor the channel is known. Three results stand out. First, after projecting out $x_0$, the GLR test declares a fault exactly when a lasso problem has a nonzero solution (a square-root lasso for unknown noise variance). Second, the null distribution of the test depends on the system only through the Gram matrix of the projected atoms, which fixes the thresholds offline. Third, for a stable plant only a bounded number of early fault times are poorly visible, whereas for undamped modes a fixed fraction of the record consists of poorly visible fault times.

\smallskip 
Future directions include faults at multiple channels, for which the covariance test \cite{LTTT14} or a detection-boundary analysis is needed; step, ramp and sensor faults, which change only the atoms; and systems with Jordan blocks.

\ifarxiv\appendices\fi

\section{Proofs}
Throughout the paper, $o(1)$ denotes a quantity that tends to $0$ as $N\to\infty$, and $o(N)$ a quantity that is $o(1)$ times $N$.

\begin{proof}[Proof of Proposition~\ref{prop:glr}]
\looseness=-1 We first compute sup over $x_0$ in the numerator and denominator of \eqref{eq:glr}.
Split $y=\Pp y+(I-\Pp)y$. Since $\Pp d_a=w_az_a$ and the two components are orthogonal, $\|y-\Ob x_0-\mu d_a\|^2=\|\Pp y-\mu w_az_a\|^2+\|(I-\Pp)(y-\mu d_a)-$ $\Ob x_0\|^2 .$ The first term does not involve $x_0$ and the second vanishes at
$x_0=\Ob^\dagger(y-\mu d_a)$, where $\Ob$ has full column rank.
This $x_0$ minimizes the exponent of the Gaussian density, and so maximizes the density,
for every $\mu$ and $\sigma^2$. Thus, 
\begin{equation}\label{eq:profilelik}
  \sup_{x_0}p(y\,|\,x_0,\mu,a,\sigma^2)=(2\pi\sigma^2)^{-Np/2}\exp\!\Big(-\tfrac{\|\Pp y-\mu w_az_a\|^2}{2\sigma^2}\Big),
\end{equation}
with $\mu=0$ under $H_0$.

\smallskip   
\looseness=-1 (i) The constant $(2\pi\sigma^2)^{-Np/2}$ in \eqref{eq:profilelik} is the same for $H_0$ and $H_1$, for known $\sigma$, so it cancels in \eqref{eq:glr}. The same holds for the density of $\Pp y$, with $r$ replacing $Np$, so $\ell=\ell_{\mathrm{proj}}$. To compute $\ell$, it we evaluate the ratio of the exponentials. By Lemma~\ref{lem:glr}, maximizing over $\mu$ turns the exponent into $-(\|\Pp y\|^2-(z_a^\top y)^2)/(2\sigma^2)$ for atom $a$ under $H_1$, and it is $-\|\Pp y\|^2/(2\sigma^2)$ under $H_0$. The log-ratio of the numerator for atom $a$ to the denominator is therefore $(z_a^\top y)^2/(2\sigma^2)$. Since the denominator does not depend on $a$ and $\log$ is increasing, $\ell=\max_a (z_a^\top y)^2/(2\sigma^2)=T/2$.

\smallskip  
(ii) (\emph{Unknown $\sigma$}) Now $\sigma^2$ is also maximized in \eqref{eq:glr}. After maximizing over $\mu$, Lemma~\ref{lem:glr} gives the exponent $-R_a/(2\sigma^2)$ under atom $a$ and $-R_0/(2\sigma^2)$ under $H_0$, with residuals $R_a=\|\Pp y\|^2-(z_a^\top y)^2$ and $R_0=\|\Pp y\|^2$. For $m>0$ and $R>0$,
\[
\sup_{\sigma^2}(2\pi\sigma^2)^{-m/2}e^{-R/(2\sigma^2)}=(2\pi eR/m)^{-m/2}.
\]
 For $\ell(y)$ we take $m=Np$, the exponent of \eqref{eq:profilelik}, so the log-ratio for atom $a$ is $\tfrac{Np}{2}\log(R_0/R_a)$. For $\ell_{\mathrm{proj}}$ we use the density of $\Pp y$, which has $m=r$, so the log-ratio is $\tfrac r2\log(R_0/R_a)$. Since $\log(R_0/R_a)=-\log\big(1-(z_a^\top y)^2/\|\Pp y\|^2\big)$ and this is increasing in $(z_a^\top y)^2$, maximizing over $a$ gives $-\log(1-U^2)$. Hence $\ell_{\mathrm{proj}}=-\tfrac r2\log(1-U^2)$ and $\ell=\tfrac{Np}{r}\,\ell_{\mathrm{proj}}$.
\end{proof}

\begin{proof}[Proof of Lemma~\ref{lem:lasso}]
\emph{Step 1: eliminate $x_0$.} Regularizers do not involve $x_0$. For fixed $\psi$, the minimum of $\|y-D\psi-\Ob x_0\|^2$ over $x_0$ equals $\|\Pp(y-D\psi)\|^2$, as in Lemma~\ref{lem:glr}. So \eqref{eq:wlasso} is equivalent to $\min_\psi\frac12\|\Pp y-\Pp D\psi\|^2+h\sum_aw_a|\psi_a|$, and \eqref{eq:sqrt} to the same problem with the norm unsquared and $\lambda$ for $h$.

\smallskip 
\looseness=-1 \emph{Step 2: normalize.} Since $w_a>0$, $\varphi_a=w_a\psi_a$ is a bijective change of variables; $\Pp D\psi=\sum_a\varphi_az_a=Z\varphi$ with $Z=[z_a]$, and the penalty becomes $h\|\varphi\|_1$. We obtain a standard lasso in the unit-norm dictionary $Z$ with data $\tilde y=\Pp y$.

\smallskip 
\looseness=-1 \emph{Step 3: optimality (P1).} The objective $\frac12\|\tilde y-Z\varphi\|^2+h\|\varphi\|_1$ is convex, so $\varphi=0$ is optimal iff $0\in-Z^\top\tilde y+h\,\partial\|\varphi\|_1|_{\varphi=0}$, that is, $\|Z^\top\tilde y\|_\infty=\max_a|z_a^\top y|\le h$, since $Z^\top\tilde y=Z^\top y$.

\smallskip 
\looseness=-1 \emph{Step 4: optimality (P2).} The objective: $\|\tilde y-Z\varphi\|+\lambda\|\varphi\|_1$. If $\tilde y\neq0$, the first norm term is differentiable at $\varphi=0$ with the gradient $-Z^\top\tilde y/\|\tilde y\|$, so, as we did in Step 3, $\varphi=0$ is optimal iff $\|Z^\top\tilde y\|_\infty/\|\tilde y\|\le\lambda$. By definition the left side is $\max_a|z_a^\top y|/\|\Pp y\|=U$, so the solution is zero iff $U\le\lambda$.
\end{proof}

\begin{proof}[Proof of (i) of Theorem~\ref{thm:long} and of Remark~\ref{rem:long}]
Write $t=\tau+1$ and $q=N-t$, so that \eqref{eq:profile} reads $1-\rho_a^2=x^\top W_N^{-1}x\,/\,b_j^\top W_qb_j$ with $x=(A^t)^\top W_qb_j$.
 Let $W_q^{1/2}$ be the symmetric square root, so $x=(A^t)^\top W_q^{1/2}\,(W_q^{1/2}b_j)$ and $b_j^\top W_qb_j=\|W_q^{1/2}b_j\|^2$. Then
\begin{align*}
&1-\rho_a^2
=\frac{\|W_N^{-1/2}(A^t)^\top W_q^{1/2}(W_q^{1/2}b_j)\|^2}{\|W_q^{1/2}b_j\|^2}\\
&\le\|W_N^{-1/2}\|^2\,\|A^t\|^2\,\|W_q^{1/2}\|^2=\frac{\lambda_{\max}(W_q)}{\lambda_{\min}(W_N)}\,\|A^t\|^2\le\kappa_\infty\,M^2\lambda^{2t}.
\end{align*}
Since $W_i$ increases with $i$ and $W_i\to W_\infty$, we have $W_q\preceq W_\infty$ and $W_N\succeq W_n\succ0$ for $N\ge n$, so $\lambda_{\max}(W_q)/\lambda_{\min}(W_N)\le\kappa_\infty$. With $\|A^t\|\le M\lambda^t$ the last line follows, which is (i).

\smallskip 
\emph{Limit and count.} For fixed $t$, $q=(N-t)\to\infty$, so $W_q,W_N\to W_\infty\succ0$, and the right side of \eqref{eq:profile} is continuous at $(W_\infty,W_\infty)$ because its denominator tends to $b_j^\top W_\infty b_j>0$. This gives $\gamma_{t,j}$. If an atom is poorly visible, $1-\rho_a^2>1-\rho_0^2$, so by (i) $t<\ln\!\big(\kappa_\infty M^2/(1-\rho_0^2)\big)/\big(2\ln(1/\lambda)\big)$, a bound free of $N$.
\end{proof}

\ifarxiv
\begin{proof}[Proof of (ii) of Theorem~\ref{thm:long}]
The idea is that, when the modes are marginal, the Gramian $W_i$ grows linearly, $W_i\approx i\bar W$, where $\bar W$ is the average Gramian per sample (defined in Step~1). The cross terms between different modes average out, and $A$ preserves the $\bar W$-norm. Then $x=(A^t)^\top W_qb_j\approx q\bar WA^{-t}b_j$, and \eqref{eq:profile} reduces to $q/N$. Write $t=\tau+1$, $q=N-t$, and let $A=\sum_\nu\nu P_\nu$ be the spectral decomposition over the distinct eigenvalues $\nu$, with $P_\nu$ the spectral projectors.

\smallskip 
\emph{Step 1: Gramian average.} For real-valued $x$, $\|CA^kx\|^2=\sum_{\nu,\nu'}(\nu\bar\nu')^k(CP_\nu x)^*(CP_{\nu'}x)$. For $\nu\neq\nu'$ on the unit circle, $\nu\bar\nu'\neq1$, so $\sum_{k<i}(\nu\bar\nu')^k$ is bounded and its average over $i$ vanishes. Hence $x^\top W_ix/i\to x^\top\bar Wx$ with $\bar W=\sum_\nu P_\nu^*C^\top CP_\nu$, which is real and symmetric as eigenvalues come in conjugates. For symmetric $M$, $e_j^\top Me_k=\tfrac14[(e_j+e_k)^\top M(e_j+e_k)-(e_j-e_k)^\top M(e_j-e_k)]$, so applying the limit to $x=e_j\pm e_k$ shows that element-wise convergence of $W_i/i$ to $\bar W$.

\smallskip 
\emph{Step 2: $\bar W$ is positive definite.} If $x^\top\bar Wx=0$, then $CP_\nu x=0$ for every $\nu$. Each $P_\nu x$ is zero or an eigenvector of $A$ lying in $\ker(C)$, which observability (PBH test) excludes. So $P_\nu x=0$ for all $\nu$ and $x=\sum_\nu P_\nu x=0$. In particular, $b_j^\top\bar Wb_j>0$.

\smallskip 
\emph{Step 3: $A$ preserves the $\bar W$-norm.} Let $M=\sup_{t\in\mathbb Z}\|A^t\|$, which is finite since $A$ is diagonalizable with unit-modulus eigenvalues. From $W_{i+1}=W_i+(A^i)^\top C^\top CA^i=C^\top C+A^\top W_iA$ we get $A^\top W_iA=W_i+(A^i)^\top C^\top CA^i-C^\top C$, and the last two terms are bounded in $i$. Dividing by $i$ gives $A^\top\bar WA=\bar W$. Hence $(A^t)^\top\bar W=\bar WA^{-t}$ and $(A^{-t})^\top\bar WA^{-t}=\bar W$ for all $t\in\mathbb Z$.

\smallskip 
\emph{Step 4: early atoms.} Fix $\delta\in(0,1)$ and consider the atoms with $q\ge\delta N$. Using Step 1, write $W_q=q\bar W+E_q$ and $W_N=N(\bar W+F_N)$ with $\|E_q\|/q\to0$ and $\|F_N\|\to0$; these errors depend on $q$ and $N$ only, not on $t$, and note that $q\ge\delta N\to\infty$. Let $s=A^{-t}b_j$, so $\|s\|\le M\|b_j\|$. By Step 3, $x=(A^t)^\top W_qb_j=q\bar Ws+e$ with $\|e\|\le M\|E_q\|\|b_j\|=o(q)$, and $s^\top\bar Ws=b_j^\top\bar Wb_j$. Putting these parts together we get, 
\begin{align*}
    x^\top W_N^{-1}x&=\frac{q^2}{N}\,b_j^\top\bar Wb_j\,(1+o(1)),\\ b_j^\top W_qb_j&=q\,b_j^\top\bar Wb_j\,(1+o(1)),
\end{align*}
uniformly in $t$ with $q\ge\delta N$, because $\|s\|$ is bounded. By \eqref{eq:profile}, $1-\rho_a^2=\tfrac qN(1+o(1))$. Since $q/N=1-t/N$, this means $|\rho_a^2-t/N|=o(1)$ uniformly over these atoms.

\emph{Step 5: late atoms.} For $q<\delta N$ we invoke the bound of part~(i), $1-\rho_a^2\le\lambda_{\max}(W_q)\|A^t\|^2/\lambda_{\min}(W_N)$. Here $\|A^t\|\le M$ and $W_q\preceq qM^2\|C\|^2I$, while $\lambda_{\min}(W_N)/N\to\lambda_{\min}(\bar W)>0$ by Step 1 and the continuity of eigenvalues. Thus $\lambda_{\min}(W_N)\ge N\lambda_{\min}(\bar W)/2$ for large $N$. So $1-\rho_a^2\le Kq/N<K\delta$ with $K=2M^4\|C\|^2/\lambda_{\min}(\bar W)$, and $|\rho_a^2-t/N|\le(1-\rho_a^2)+q/N\le(K+1)\delta$.

\smallskip 
Combining Steps 4 and 5, $\limsup_N\max_\tau|\rho_a^2-t/N|\le(K+1)\delta$ for every $\delta$, which proves \eqref{eq:marginal}. The count follows as for large $N$, atoms with $t/N<\rho_0^2-\epsilon$ are poorly visible and those with $t/N>\rho_0^2+\epsilon$ are not, so their count is $\rho_0^2N+o(N)$.
\end{proof}
\fi

\ifarxiv
\section{Additional Simulations}\label{app:sims}
\emph{(Visibility versus horizon)}
Table~\ref{tab:scaling} counts the poorly visible atoms, $\#\{\tau:\rho_\tau^2<0.9\}$, as the record length $N$ grows. For strictly stable systems the count does not depend on $N$. It matches the window $\ln(10)/(2\ln(1/|a_0|))$ of the scalar case in Section~\ref{sub:profile}, which is $10.9$ for $A=0.9$ ($10$ counted) and $114.6$ for $A=0.99$ ($114$ counted at $N=400$). For $N\le200$ the record is shorter than this window, so the count is smaller.

For systems with unit-modulus modes the count grows like $0.9N$, as Theorem~\ref{thm:long}(ii) predicts. For the integrator, $\rho_\tau^2=(\tau+1)/N$, so the count is the number of $\tau$ with $(\tau+1)/N<0.9$, namely $89$, $179$ and $359$ for $N=100$, $200$ and $400$. A four-state, two-output system with two planar rotations, by $0.7$ and $1.9$ radians, gives nearly the same counts. For the rotation with $|\lambda|=1$, the largest gap $|\rho_\tau^2-(\tau+1)/N|$ over the atoms with $N-\tau-1\ge0.1N$ is $0.015$, $0.0074$, $0.0037$ and $0.0019$ for $N=100$, $200$, $400$ and $800$. The gap halves each time $N$ doubles, so it decays like $1/N$.

For part (i) we use a random stable four-state system with $p=2$ outputs and spectral radius $0.92$. Formula~\eqref{eq:profile} matches the direct projection to $7\times10^{-16}$. The bound of Theorem~\ref{thm:long}(i) holds, and it is loose: $1-\rho_a^2$ stays below $\kappa(W_N)\|A^t\|^2$, and the largest ratio of the two is $1.7\times10^{-3}$.

\begin{table}[t]
\centering\small\setlength{\tabcolsep}{6pt}\caption{Poorly visible atoms $\#\{\tau:\rho_\tau^2<0.9\}$.}\label{tab:scaling}
\begin{tabular}{lccc}\toprule
system & $N{=}100$ & $N{=}200$ & $N{=}400$\\\midrule
$A=1$ (integrator) & 89 & 179 & 359\\
$A=0.99$ & 75 & 107 & 114\\
$A=0.9$ & 10 & 10 & 10\\
rotation, $|\lambda|=1$ & 91 & 180 & 360\\
rotation, $|\lambda|=0.95$ & 23 & 23 & 23\\\bottomrule
\end{tabular}
\end{table}

\emph{({Thresholds})} Table~\ref{tab:calib} gives the false-alarm probability of the max test at nominal value $0.05$ ($10^5$ runs, $\sigma=1$). The Rice bound is not uniformly tighter than Bonferroni, and the independence approximation of Fuchs (\v{S}id\'ak) coincides with Bonferroni to within $0.001$.

\begin{table}[t]
\centering\small\setlength{\tabcolsep}{6pt}\caption{False-alarm probability of the max test at nominal $0.05$.}\label{tab:calib}
\begin{tabular}{lcccc}\toprule
system & Bonferroni & \v{S}id\'ak & Rice & exact\\\midrule
$A=1$ & 0.014 & 0.014 & 0.031 & 0.049\\
$A=0.9$ & 0.028 & 0.029 & 0.039 & 0.051\\
rotation, $|\lambda|{=}0.9$ & 0.041 & 0.041 & 0.034 & 0.051\\\bottomrule
\end{tabular}
\end{table}

\emph{(Two inputs and two outputs)}
Table~\ref{tab:mimo} uses the system of Section~\ref{sec:sims} with $\beta=0$, $0.3$ and $1$, for which the spectral radius is $\rho(A)=1$, $0.972$ and $0.909$, respectively. For $\beta>0$ the number of poorly visible atoms does not depend on $N$. For $\beta=0$ it grows like $0.9N$. For $\beta=0.3$ and $\beta=1$, the window of Section~\ref{sub:profile} computed from the slowest mode is $40.2$ and $12.1$ samples, against $37$ and $13$ atoms counted. The $U$ test, run with unknown $\sigma=3$ and $x_0\sim\mathcal N(0,20^2I_4)$, has false-alarm probability between $0.048$ and $0.053$. The Bonferroni threshold gives $0.014$ to $0.018$, and the \v{S}id\'ak threshold agrees with it to within $0.001$. Both thresholds are conservative.

\begin{table}[t]
\centering\small\setlength{\tabcolsep}{6pt}
\caption{Poorly visible atoms on channel~1 (top) and false-alarm probability at nominal value $0.05$ (bottom).}\label{tab:mimo}\label{tab:mimo}
\begin{tabular}{lccc}\toprule
 & $\beta=0$ & $\beta=0.3$ & $\beta=1$\\\midrule
$N=100$ & 90 & 37 & 13\\
$N=200$ & 179 & 37 & 13\\
$N=400$ & 360 & 37 & 13\\\midrule
$U$ test & 0.048 & 0.053 & 0.053\\
Bonferroni & 0.017 & 0.014 & 0.018\\\bottomrule
\end{tabular}
\end{table}

\begin{figure}[t]
\centering\includegraphics[width=\columnwidth]{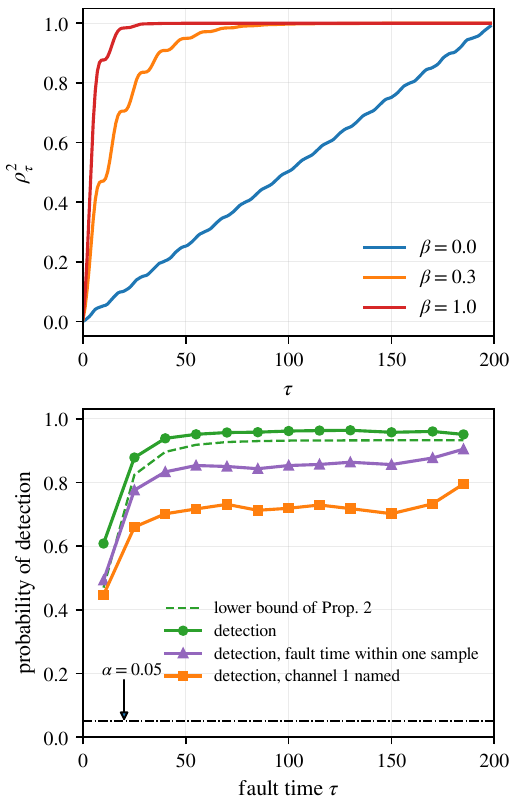}
\caption{Two-mass chain with damping $\beta K$, $N=200$, fault on channel~1 at time $\tau$. Top: visible fraction $\rho_\tau^2$ against $\tau$ for $\beta=0$, $0.3$ and $1$, respectively. Bottom: for $\beta=0.3$, $\mu\|d_a\|/\sigma=5$ and $3000$ trials per $\tau$, the probability that the max test detects the fault (green, lower bound of Proposition~\ref{prop:pd} dashed), that it detects and estimates the fault time within one sample (purple), and that it detects and names channel~1 (orange). The dash-dotted line is $\alpha$.}\label{fig:mimofull}
\end{figure}
\fi
\section*{Acknowledgment}
The author used a generative AI tool (Claude, Anthropic) for writing and coding assistance. The author is responsible for all results and claims in the paper.
\bibliographystyle{IEEEtran}
\bibliography{references}
\end{document}